\documentclass[onefignum,onetabnum]{siamart190516}

\usepackage{url}

\def\R{\bf R}

\usepackage{lipsum}
\usepackage{amsfonts}
\usepackage{graphicx}
\usepackage{epstopdf}
\usepackage{algorithmic}
\ifpdf
  \DeclareGraphicsExtensions{.eps,.pdf,.png,.jpg}
\else
  \DeclareGraphicsExtensions{.eps}
\fi

\newsiamremark{remark}{Remark}
\newsiamremark{hypothesis}{Hypothesis}
\crefname{hypothesis}{Hypothesis}{Hypotheses}
\newsiamthm{claim}{Claim}

\headers{Granular Sieving for Deterministic Global Optimization}{T. Qian, L. Dai, L.Zhang, and Z. Chen}

\title{A Granular Sieving Algorithm for Deterministic Global Optimization
\thanks{Submitted to the editors DATE. The first author invented the theory. The second author,being considered to give equal contribution to the paper, realized the algorithm and conducted the enormous experimental examples. 
\funding{This work was funded by the Science and Technology Development Fund of Macau SAR 0123/2018/A3 and Multi-year Research Grant MYRG2018-00111-FST.}}}

\author{Tao Qian\thanks{Macao Center of Mathematical Sciences,
		Macau University of Science and Technology,
		Taipa,
		Macau,
		China 
  (\email{tqian@must.edu.mo}).}
\and Lei Dai\thanks{Faculty of Science and Technology,
	University of Macau,
	Taipa,
	Macau,
	China 
  (\email{yb87427@um.edu.mo}, \email{lmzhang@um.edu.mo}).}
\and Liming Zhang\footnotemark[3] \thanks{Corresponding author}
\and Zehua Chen\thanks{College of Data Science,
Taiyuan University of Technology,
Taiyuan,
Shanxi,
China
(\email{chenzehua@tyut.edu.cn}).}}

\usepackage{amsopn}

\ifpdf
\hypersetup{
	pdftitle={A Granular Sieving Algorithm for Deterministic Global Optimization},
	pdfauthor={Tao Qian, Lei Dai, Liming Zhang, and Zehua Chen}
}
\fi

\begin{document}
	
	\maketitle
	
	\begin{abstract}
		A gradient-free deterministic method is developed to solve global optimization problems for Lipschitz continuous functions defined in arbitrary path-wise connected compact sets in Euclidean spaces. The method can be regarded as granular sieving with synchronous analysis in both the domain and range of the objective function. With straightforward mathematical formulation applicable to both univariate and multivariate objective functions, the global minimum value and all the global minimizers are located through two decreasing sequences of compact sets in, respectively, the domain and range spaces. The algorithm is easy to implement with moderate computational cost. The method is tested against extensive benchmark functions in the literature. The experimental results show remarkable effectiveness and applicability of the algorithm.
	\end{abstract}
	
	
	
	\section{Introduction}
	\label{intro}
	 In the literature, many theoretical and practical scientific and engineering problems can be reduced, or are closely related, to various types of optimization problems \cite{Floudas2005,horst2013,pinter2006,Torn1989}. The task of global optimization is to find the global minimum or maximum (conventionally, the minimum) of an objective function in some region of interest \cite{floudas2009,Horst2000,pardalos2000}. Global optimization algorithms are generally divided into the deterministic and the stochastic types \cite{horst2013,liberti2005}. Deterministic methods theoretically guarantee that at least one or all global solutions (minimizers) can be found \cite{barhen1997,Floudas2013,Horst20132,niebling2019}, while stochastic methods find  the solutions in the sense of probability \cite{jones1998,liberti2005,zhigljavsky2007,Nemirovski2009}.
	In the literature, remarkable theoretical and practical deterministic methods have been developed \cite{booker1999,huyer2008,jones2001,jones1993,pardalos2000,powell1999,Rios2013}. Stochastic methods take random selections of the initial values and do probability analysis on the randomization results \cite{horst2013,liberti2005,jones1998,zhigljavsky2007}. For many optimization methods the initial value to start the algorithm is crucial. An algorithm involving extensively many initial value selections leads, in fact, a computational NP hard problem. In such a sense, stochastic methods are manipulation and amendments of the un-solvability of the initial value NP hard problems, offering acceptable accuracy and efficiency \cite{Floudas2013,liberti2005,Nesterov2012}. 
	
	In this study, we introduce a practical and efficient deterministic method for general global optimization. As a well designed  exhaustive search method, it is applicable to objective functions of all dimensions and theoretically ensures finding, with fast convergence and satisfactory precision, of the global minimum value together with the entire set of the global minimizers. The proposed algorithm is inspired partially by our previous study on a particular type of non-Lipschitzian functions but with theoretically accurate solutions for the optimal global $n$-parameter \cite{Qian2014,Qian2019,Qian2013}, and partially by the granular computing idea in artificial intelligence \cite{bargiela2016,Lin2009,Pedrycz2008}. The algorithm is thus named granular sieving (GrS). Through enhancing the primary granular computation idea, the proposed GrS algorithm yields a complete solution set. The proof of validity of the algorithm is fundamental and explicit that can be seen as a refinement version of existence and attainability of the global minimum value of a continuous function \cite{jones1993,Pinter2013} defined on a compact (bounded and closed) set in the standard analysis textbooks.
	
	\def\bR{\bf R}
	\def\bx{\bf x}
	The objective functions that we are aiming to study are of the form
	$f(x_1,\cdots,x_n),$ where $f$ is a real-valued continuous function with the variables $x_1,\cdots,x_n$ ranging in a compact (bonded and closed) set $C$ of a real- or even complex- $n$-dimensional Euclidean space. Examples of such functions are everywhere, including distance from a single point to a domain in a metric space, energy functions, minimum surface problems, costing functions with parameters in machine learning practices, etc. In most existing algorithms local comparison is employed. In one way or another a local comparison method involves gradients. The GrS algorithm we propose assumes that the Lipschitz condition is satisfied, which is weaker than existence and boundedness of the gradients of the objective function. The Lipschitz condition stands for the following: There
	exists a constant ${M}$ such that
	\begin{eqnarray}\label{Lip}
		\frac{|f(x+\Delta x)-f(x)|}{|\Delta x|}\le {M}. \end{eqnarray}
	We call the difference-quotient quantity on the left-hand-side of 
	(\ref{Lip}) as \emph{relative variation}.
	We note that the GrS algorithm we are proposing is implementable and effective for all upper bounds $
	{M}$ defined by Eq. (\ref{Lip}), but not restricted to any specific one such as the infimum of the upper bounds $M$. In the sequel the infimum is denoted by $M_{0}$ and called the Lipschitz constant. It can be shown that
	\begin{eqnarray*}\label{M0}
		M_0&=&\min\{ M\ :\ M\ \rm is \ an \ upper\ bound\ of \ the\ relative \ variations \}\\
		&=&\sup\{\frac{|f(x+\Delta x)-f(x)|}{|\Delta x|}\ :\ x,\ x+\Delta x\in C, \Delta x\ne 0\}.
	\end{eqnarray*}
	A point $x^{\ast}$ is said to be a \emph{global minimizer} of the function $f(x)$ defined on $C,$ if $x^{\ast}\in C$ and $f(x^{\ast})\le f(x)$ for all $x\in C.$ In such case $f(x^\ast)$ is called the \emph{global minimum value} of $f$ on $C$ and denoted as $f_{\min}.$ In below we denote $v=f_{\min}.$
	
	It would be of great convenience if one could find and hence use the
	precise bound $M_0,$ that is the least upper bound $M$ of the relative variations. Finding the least upper bound of the relative variations would often be a problem of similar difficulty to the original optimization problem. The algorithm that we propose, however, does not require precise knowledge of the least upper bound $M_0$ in question.
	
	In proving convergence of the algorithm, the compactness assumption of the domain of the objective function $f$ is necessary. It is seen, however, optimization problems of objective functions defined on, say, open or unbounded sets, can possibly also employ the proposed algorithm, if the quantity of their boundary values are well controlled.
	
	The method is tested against an extensive number of benchmark functions and, in particular, compared with some 12 popular methods on the top 10 hardest functions listed in the test function index \cite{testfunctionindex}. The results show remarkable effectiveness and applicability. With many of the tested famous and typical examples, GrS surpasses the popular stochastic algorithms. We anticipate that the GrS method would be applicable in various types of theoretical and practical problems and arouse certain attention and interest.

	\section{The GrS Algorithm Outline}
	We will be working with a real-valued continuous objective function with real variables. For complex variables the theory and algorithm are similar.
	Let $C,$ being the domain of the objective function under study, be a compact (closed and bounded) path-wise connected set in the Euclidean space $\mathbb{R}^n.$
	The algorithm will construct a decreasing sequence of compact sets $C=C_1\supset \cdots \supset C_k\supset \cdots$ in
	the domain $C$ of the objective function with the following properties.
	For each $k=1,\cdots,$ one constructs a partition of $C_k,$
	$\mathbb{N}_k=\{C^{(k)}_l\ :\ l=1,\cdots,l_k\},$ such that
	$C_k=\cup_{l=1}^{l_k}C^{(k)}_l.$
	For each fixed $k$ the pairs $(k,l), l=1,\cdots,l_k,$
	and no others, are called $k$-\emph{admissible sets}, or just \emph{admissible sets}, in brief.
	The sets $C^{(k)}_l$ satisfy
	the following conditions:\\
	
	(i) When $k$ grows larger, $\mathbb{N}_k$ becomes finer in the sense that a set $C^{(k+1)}_{l'}$
	is either entirely contained in $C^{(k)}_{l}$ or has empty overlap with the interior of $C^{(k)}_{l}.$\\
	
	(ii) The maximal diameter of $C^{(k)}_l, l=1,\cdots,l_k,$ denoted as $\delta_k,$
	tends to zero along with
	$k\to \infty.$ The diameter is defined as the maximal distance between any two points in the set.\\
	
	(iii) Each set $C^{(k)}_l$ is compact whose boundary
	can have non-empty overlap with the boundaries of the other $C^{(k')}_{l'}$. \\
	
	(iv) Along with the construction we name in each $C^{(k)}_l$ a representative point $x_{l}^{(k)}$.  When $C^{(k)}_l$ is geometrically symmetric we may select $x_{l}^{(k)}$ as the geometrical center of $C^{(k)}_l.$ For this reason we, in below, call $x_{l}^{(k)}$ as $\lq\lq$center" of $C^{(k)}_l,$ although they are not necessarily the centers, but only representative points of $C^{(k)}_l.$\\
	
	We will be using mappings between sets defined as
	
	\[ f(A)=\{b\in {\R} \ :\ \exists a\in A, b=f(a)\}, \quad A\subset C,\]
	where ${\R}$ stands for the set of real numbers.
	
	As a fundamental result of a continuous function defined in a compact set,
	there exists a non-empty compact set $\tilde{C}$ on
	which the function takes the global minimum
	value $v,$ that is $f(\tilde{C})=\{v\}.$ The purpose of the algorithm is to inductively construct the set $C_{k+1}$ through some sets in the partition $\mathbb{N}_{k}=\{C^{(k)}_l\ :\ l=1,\cdots,l_k\}$ of $C_k,$ and to realize $\tilde{C}=\cap_{k=1}^\infty C_k,$ and $f(\tilde{C})=\{v\}.$  \\
	
\begin{figure}[h]
		\centering
		\includegraphics[scale=0.25]{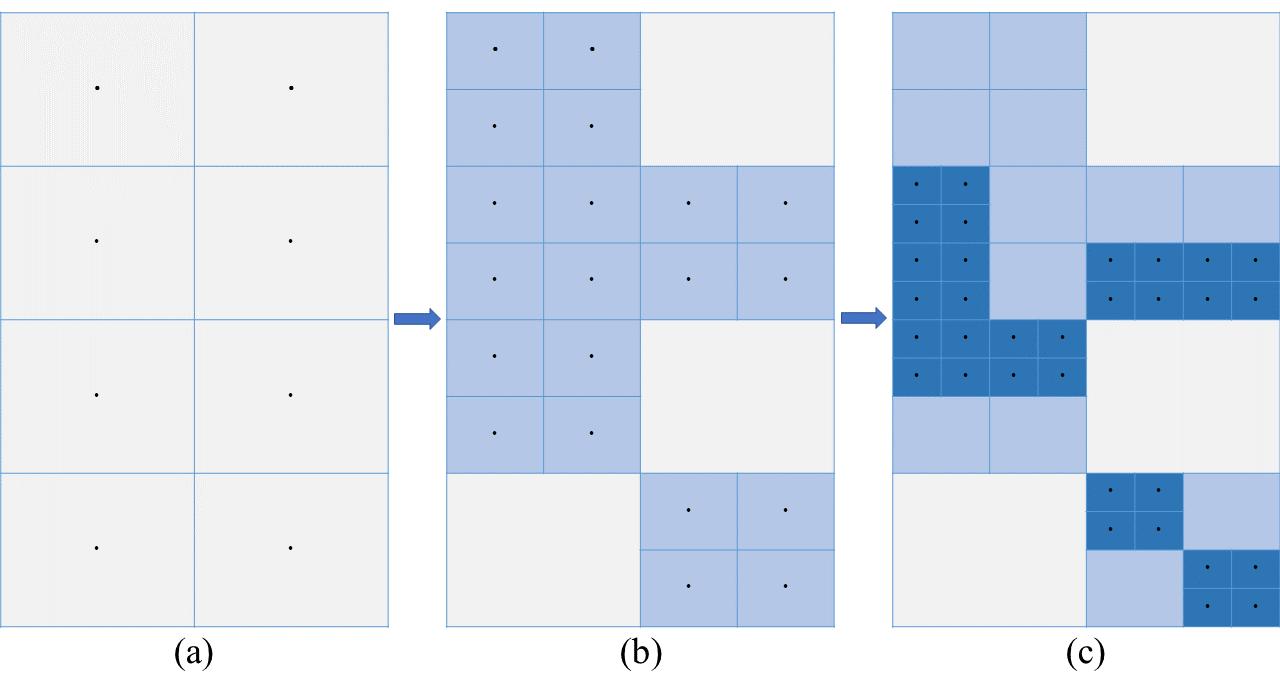}
		\caption{An example of the first three consecutive partition processes of a two-dimensional function. The black dots represent the centers of the boxes.}
		\label{f1}
	\end{figure}

	\noindent{\bf The Algorithm Details} The proposed GrS algorithm for finding
	the global minimum value $v$ and at the same time all the global minimizers is of a fundamental nature, as given below. For the global maximum the algorithm and the theory are the same. The idea of the algorithm is as follows: For each fixed level $k$ the set $C_k$ is expressible by the sets in its partition: $C_k=\cup_{l=1}^{l_k}C^{(k)}_l, C^{(k)}_l\in \mathbb{N}_{k}, l=1,\cdots,l_k.$ For the fixed $k,$ exam all the function values at the center points $x^{(k)}_l$ of the $k$-level sets $C^{(k)}_l,$ and select the minimum one among the function values $f(x^{(k)}_l)$, and denote it by $v_k.$ If a center $x^{(k)}_{l'}$ satisfies the relation $f(x^{(k)}_{l'})>v_k,$ we can determine by using the Lipschitz condition whether there would be some function values $f(x), x\in C^{(k)}_{l'},$ that can possibly be equal to, or below $v_k.$ The $k$-level sets $C^{(k)}_{l'}$ that are of such possibility will be called $\lq\lq$good sets". If there are no function values that could be equal to or below the value $v_k,$ then we simply delete the whole set $C^{(k)}_{l'}.$ The deleted sets $C^{(k)}_{l'}$ are referred as $\lq\lq$bad sets". The union of the $\lq\lq$good sets" is defined to be our next generation set $C_{k+1}.$ The algorithm consists of iteration of the above sorting process that forms $C_{k+1}$ as union of the good subsets in $\mathbb{N}_{k}$ of the partition of $C_k.$ Figure \ref{f1} shows, as an example, the first three consecutive partitions of a two-dimensional function. Figure \ref{f1}(a) shows the first partition, and Figure \ref{f1}(b) and \ref{f1}(c) illustrate that only those classified as $\lq\lq$good sets" are collected and further partitioned, expressed as the light-blue boxes in Figure \ref{f1}(b) and the dark-blue boxes in Figure \ref{f1}(c). Finally one takes
	$\cap_{k=1}^\infty C_k,$ which is non-empty giving rise to the global minimum $v.$ Below is a more detailed explanation of the algorithm.
	
	Carrying out the above mentioned sorting at the level $k=1$ is based on constructing
	a partition $\mathbb{N}_1$ of $C=C_1.$  The partition should be fine enough so that the maximal diameter $\delta_1$ of all $C^{(1)}_l$ makes the quantity $\delta_1M_0$ strictly smaller than the largest possible variation of the function values $f(x^{(1)}_l)$ from their minimum value $v_1,$ where
	\begin{eqnarray}\label{v1}
		v_1=\min \{f(x^{(1)}_l)\ :\ l=1,\cdots,l_1\}.
	\end{eqnarray}
	We divide the sets in the $C_1$-partition $\mathbb{N}_1$ into two non-overlap subclasses:
	$\mathbb{N}_1=\mathbb{B}_1\cup \mathbb{G}_1$. \\
	
	The subclass
	$\mathbb{B}_1,$ as the collection of the $\lq\lq$bad" sets, consists of the sets
	$C^{(1)}_l$ whose centers $x^{(1)}_l$ satisfy
	\begin{eqnarray}\label{first1}
		f(x^{(1)}_l)> v_1+\delta_1M_0.
	\end{eqnarray}
	
	Notice that for any $C^{(1)}_l$ the quantity $\delta_1M_0$ is the greatest possible variation of the function values in the set $C^{(1)}_l$ from its center.  The inequality (\ref{first1}) sorts out those $C^{(1)}_l$ in which any function value cannot reach
	the minimum level $v_1$, let alone below it. The corresponding set $C^{(1)}_l$ then should be deleted. If no sets are classified into $\mathbb{B}_1,$ it means that the partition $\mathbb{N}_{1}$ is not fine enough, or $\delta_1$ is not fine enough.\\
	
	The subclass $\mathbb{G}_1,$ as the collection of the $\lq\lq$good" sets,
	consists
	of those $C^{(1)}_l$  whose centers $x^{(1)}_l$ satisfy
	the opposite inequality
	\begin{eqnarray}\label{first2}
		f(x^{(1)}_l)\leq v_1+\delta_1 M_0.
	\end{eqnarray}
	If a set $C^{(1)}_l$ belongs to
	the collection $\mathbb{G}_1,$ due to the inequality (\ref{first2}), one cannot exclude the possibility that some function values of $x\in C^{(1)}_l$ can possibly reach or even below the minimum $v_k.$ Those sets $C^{(1)}_l$ should be kept and further analyzed.\\
	
	We eliminate all the sets in the subclass $\mathbb{B}_1$ of the bad sets  and
	union those in the subclass $\mathbb{G}_1$ of  the good sets.  Define $C_2=\cup \{C^{(1)}_l\ :\ C^{(1)}_l\in \mathbb{G}_1\},$ and construct a partition
	$\mathbb{N}_2$ of $C_2$ finer than what is inherited from $\mathbb{N}_1$.\\
	
\begin{figure}[h]
		\centering
		\includegraphics[scale=0.50]{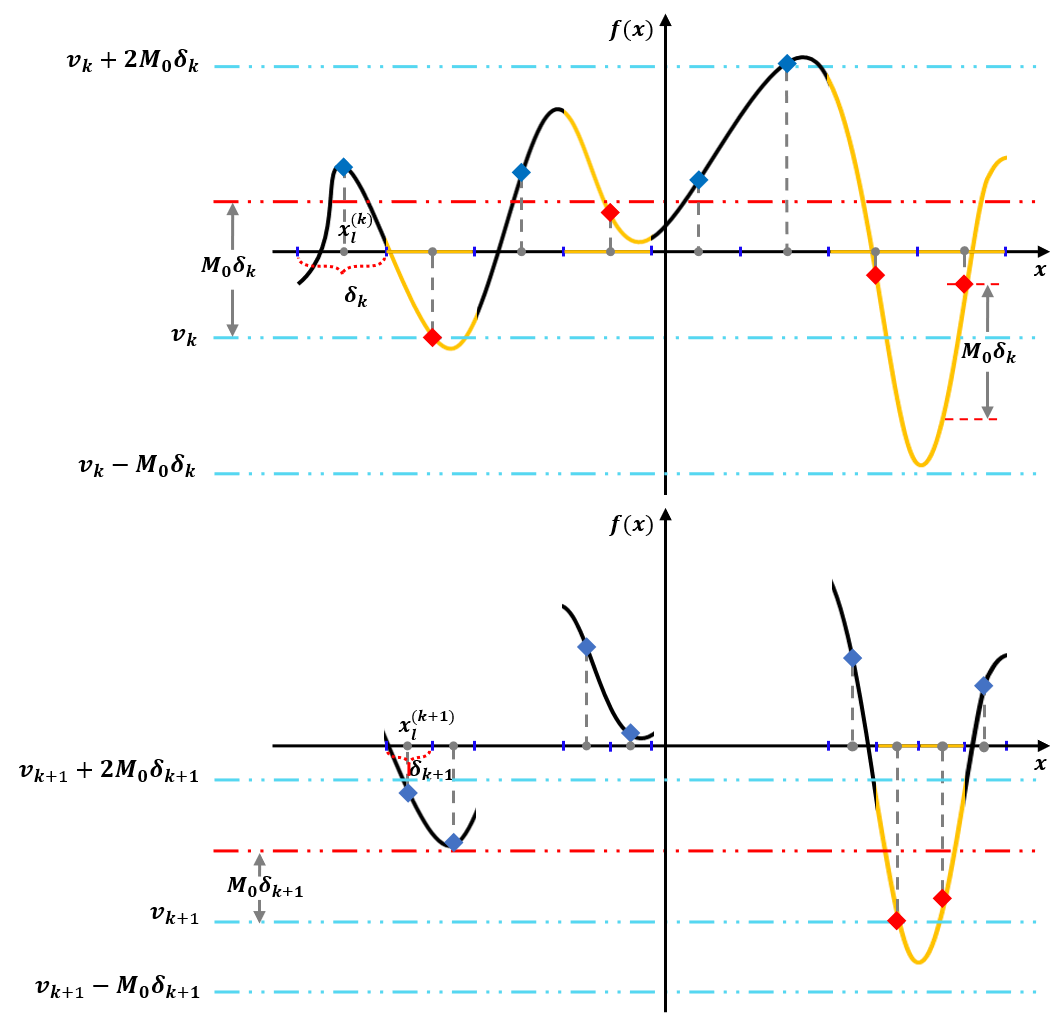}
		\caption{Demonstration diagram of the granular-sieving process of a one-dimensional function at the $k$-th and $(k+1)$-th partition levels. On the x axis, the gray dots are centers of the granule sets in the partition sets. The yellow and black intervals and function graphs represent, respectively, the $\lq\lq$good" and $\lq\lq$bad" sets in the domain and the corresponding function values. The red and blue diamonds represent the function values of the partition centers of the $\lq\lq$good" and $\lq\lq$bad" sets, respectively. The $\lq\lq$good sets" are those with partition centers below the red dotted line, while the $\lq\lq$bad sets" are those with partition centers above the red dotted line.}
		\label{f2}
	\end{figure}

	With a suitable partition $\mathbb{N}_2, C_2=\cup_{l=1}^{l_2}C^{(2)}_l.$ Let $\delta_2$ be the maximal diameter of the subsets in the partition $\mathbb{N}_2$. We define
	\begin{eqnarray}\label{v2}
		v_2=\min \{f(x^{(2)}_l)\ :\ l=1,\cdots,l_2\}.
	\end{eqnarray}
	Using the criteria 
	\begin{eqnarray}\label{second1}
		f(x^{(2)}_l)> v_2+\delta_2 M_0
	\end{eqnarray}
	and
	\begin{eqnarray}
		\label{second2} f(x^{(2)}_l) \leq v_2+\delta_2 M_0
	\end{eqnarray}
	we divide the sets $C^{(2)}_l$ in the partition $\mathbb{N}_2$  into a bad set subclass $\mathbb{B}_2$ and a good set subclass $\mathbb{G}_2,$ respectively.
	Eliminate all the bad sets $C^{(2)}_l$ whose centers satisfying the inequality
	(\ref{second1}) and keep all the good sets as the collection $\mathbb{G}_2.$ Define
	$C_3=\cup \{C^{(2)}_l\ :\ C^{(2)}_l\in \mathbb{G}_2\}$ and construct
	$\mathbb{N}_3$ as a finer partition of $C_3,$ and so on. Figure \ref{f2} shows the granular sieving process of a one-dimensional function at the $k$-th and$(k+1)$-th partition levels.

	Repeating this process inductively with 
	$C_{k+1}=\cup_{C_{l}^{(k)}\in \mathbb{G}_k}C^{(k)}_l$ and making use the condition $\delta_k\to 0,$ one has
	
	\begin{theorem}\label{Th1}
		$\cap_{k=1}^\infty C_k=\tilde{C}\ne \emptyset,$ and $f(\tilde{C})=\{v\}.$
	\end{theorem}
	
	\noindent{\bf Proof of Theorem \ref{Th1}}  We note that for each $k,$ the set $C_{l}^{(k)}\in \mathbb{G}_k$ if and only if 
	\begin{eqnarray}\label{3}
		f(x_{l}^{(k)})\le v_k+\delta_{k}M_{0}.
	\end{eqnarray}
	We also have,
	\begin{eqnarray}\label{4}
		|f(x)-f(x_{l}^{(k)})|\leq \delta_{k}M_{0}.
	\end{eqnarray}
	They imply
	\[ f(C_{l}^{(k)})\subset [v_k-\delta_kM_0,v_k+2\delta_kM_0]=I_k^{M_0}.\]
	Therefore, 
	\[ f(C_{k+1})=f(\cup_{C_{l}^{(k)}\in \mathbb{G}_k}C_{l}^{(k)})=\cup_{l=1}^{l_k}f(C_{l}^{(k)})\subset
	I_k^{M_0}.\]
	There holds
	\[ f(\cap_{k=1}^\infty C_{k+1})=\cap_{k=1}^\infty f(C_{k+1})\subset \cap_{k=1}^\infty I_k^{M_0}.\]
	Since the compact sets sequence $C_2\supset C_3\supset\cdots C_{k+1}\supset\cdots$ 
	has the finite-intersection-non-empty property, $C\supset\tilde{C}=\cap_{k=1}^\infty C_{k+1}\ne \emptyset$ \cite{hewitt2013}. On the other hand, $v\in 
	I_k^{M_0}$ for each $k,$ and $|I_k^{M_0}|=3\delta_kM_0\to 0,$ we have $\{v\}=\cap_{k=1}^\infty I_k^{M_0}.$ Therefore,
	\[ f(\tilde{C})=\{v\}.\]
	The proof is complete. \\

	\section{Further Algorithm Analysis}
	\label{sec:3}

	\noindent{\bf Degree of Finesse of the Partitions }
	The above process requires the partitions $\mathbb{N}_1, \mathbb{N}_2,\cdots ,\mathbb{N}_k,\cdots $ to satisfy that for each $k$ both the inequalities
	\begin{eqnarray}
		f(x^{(k)}_l)-v_k\leq \delta_k M_0 \quad {\rm and}\quad
		f(x^{(k)}_l)-v_k>\delta_k M_0
	\end{eqnarray}
	have solutions for some $l=1,\cdots,l_k.$  Analytically, this requires $\delta_k$ to satisfy
	\begin{eqnarray}\label{ine}
		A_k<\delta_k<B_k,
	\end{eqnarray}
	where
	\[ A_k=\frac{\min\{f(x^{(k)}_l)-v_k : l=1,\cdots,l_k\}}{M_0}\]
	and
	\[ B_k=\frac{\max\{f(x^{(k)}_l)-v_k : l=1,\cdots,l_k\}}{M_0}.\]\\
	
	When $\delta_k$ satisfies (\ref{ine}) and is small so to get close to $A_k$ one has more bad sets to delete and less good sets to keep. In such case in the $k$-iteration step the positions of the global minimizers are more accurately located, but as compensation the algorithm has a greater computational cost. When $\delta_k$ is large so to get close to $B_k$ one has more good sets to keep and less bad sets to delete. In such case at the $k$-iteration step the global minimizers are not so well located but the computation cost is lower.\\

	As analyzed in the last section the design of the partition $\mathbb{N}_k$ is crucial. A design needs to balance the computation complexity, the algorithm effectiveness and its efficiency: At each sorting process one wishes to delete a sufficient amount of bad sets so to quickly reduce the area of the set containing the minimizers. In practice, for simplicity we divide the remaining set $C_k,  (k>1)$ into the same number of equal parts and take $\delta_k$ to be one proportional to $m^{-k}$ for some integer $m$ (see the algorithm design in the following sections). \\
	
	\noindent{\bf Pseudo-Upper Bounds of the Algorithm}
	An important issue to be addressed is related to the values $M$ that are actually used in the algorithm.
	In the algorithm the least upper bound $M_0$ of the bounds of the relative variations is crucial but maybe unknown in the practice. 
	Determination
	of the least bound $M_0,$ or even any bound $M$ is a similar problem.
	This is, in fact, a common problem for all Lipschitz-constant-dependent algorithms: If the Lipschitz constant is not known, one does not know what $M$ should be used in the algorithm.\\
	
	This $M_0$-paradox may be solved through a concept called  \emph{empirical}-, or \emph{pseudo}-$M$, that we now formulate.  We can use an arbitrary positive number $M>0,$ not necessarily a bound of the relative variations, to replace $M_{0}$ in
	the algorithm described by Theorem 1. We call
	such attempted $M$ a \emph{pseudo}-$M$. We will show that the algorithm with a pseudo-$M,$ as substitution of $M_0,$ may still be carried on and convergent. The corresponding algorithm is denoted as $A_{M}$ that converges on the domain of $f$ side to a set $\tilde{C}_{M},$ called the \emph{set of $M$-pseudo minimizers of} $f;$ and a single point $v^M$ on the range of $f$ side, called \emph{the $M$-pseudo minimum of} $f.$
	The
	following theorem reveals the dynamic relation between the pseudo-$M$'s
	and their corresponding algorithm outcomes.\\
	
	\begin{theorem}\label{Th2} (i) For any $M>0$, the algorithm $A_{M}$
		converges to a compact set $\tilde{C}^{M}$ and a
		single value $v^{M}$ such that $f(\tilde{C}^{M})=\{v^M\};$ and, (ii) If $M_{1}<M_{2},$ then $v^{M_1}\ge v^{M_2}.$ Consequently, if $M\ge M_{0},$
		then $\tilde{C}^{M}=\tilde{C}^{M_0}=\tilde{C},$ and $f(\tilde{C}^{M})=\{v^{M_0}\}=\{v\}.$
	\end{theorem}
	
	\noindent{\bf Proof of Theorem \ref{Th2}} (i).
	No matter how large or small $M$ is, the algorithm procedure given in Theorem \ref{Th1} converges on the domain side to a non-empty compact set $\tilde{C}^M,$ and on the range side to a single point $v^M.$ In fact, in the self-explanatory notation, due to the non-empty intersection property of the sequence of the compact sets $\{C_k^M\}_{k=1}^\infty,$ there holds $\emptyset\ne \tilde{C}^M=\cap_{k=1}^\infty C_k^M,$ and
	\[ \emptyset\ne f(\tilde{C}^M)=f(\cap_{k=1}^\infty C_k^M)\subset \cap_{k=1}^\infty f(C_k^M)\subset \cap_{k=1}^\infty I_k^M.\]
	The set $f(\tilde{C}^M)$ can contain only one point, $v^M,$ because 
	\[ \lim_{k\to \infty}|I_k^M|=\lim_{k\to \infty}3\delta_kM=0.\]
	(ii). First we note that the results of $A_M$ is independent of the partitions in use. Assume $M_1<M_2.$ By using the same partitions at every $k$-step, with the pseudo-bounds $M_1, M_2,$ respectively, we throw away the bad sets $C^{M_j,(k)}_l$ satisfying the condition
	\[f(x_{l}^{(k)})> v_k+\delta_{k}M_{j}, \qquad j=1,2.\]
	Since for $j=1$ the set $C^{M_j,(k)}_l$ is easier to satisfy the above inequality, or, in other words, easier to be bad, hence $C^{M_1}_{k+1}$ collects less sets than $C^{M_2}_{k+1}.$ As a consequence, $v_{k+1}^{M_1}\ge v_{k+1}^{M_2}.$ By taking limit $k\to \infty,$ we have $v^{M_1}\ge v^{M_2}.$ The proof is complete.\\
	
	On one hand, $M_0$ is difficult to find. On the other hand, 
	$M_0$ is not necessarily to be found. 
	In order to
	produce the correct global minimum value of the objective function a pseudo-$M$ does not necessarily reach $M_{0}.$
	It is easily observed that large relative variations may not present
	around the global, or even local, minimizers or maximizers. 
	
	To treat the paradox, the GrS algorithm is combined with a pseudo-$M$ scheme involving a finite sequence of pseudo-$M_i$'s, that is $M_1<\dots<M_n$. The starting $M_1$ and the ending $M_n,$ however, are the art part of the method. $M_1$ in our examples may be estimated through the first partitions. With concrete problems although gradients are not directly used, their global bounds may be used to derive a practical bound of the relative variations. Such example include $n$-best kernel approximation in reproducing kernel Hilbert spaces (\cite{saitoh2016} and references therein where the open domain of the objective function can also be reduced to a compact domain set).  The stopping time $M_n$ may be decided based on observation of non-improving of the optimizing result: when the values $v^{M_n}$ become stable one can decide to end the process. Actual use of pseudo-$M$'s in concrete questions is empirical, crucial, and is the true art part of the algorithm that requires further and deep studies.
	
	In practice, solving a global optimization problem is a trade-off problem: It is to find balance between accuracy and cost \cite{barhen1997}. When employing the same set of consecutive partitions for a larger $M$, the convergence of $\mathcal{A}_M$ may be slower due to the fact that at each iteration step
	a less number of bad sets are thrown away. Larger pseudo-$M$'s result in better accuracy but rise higher computational costs. Small pseudo-$M$'s on the other hand may not give the right results but have fast convergence.
	
	\section{Experimental results}

	\noindent{\bf Test functions}
	Typical test functions for general optimization problems are provided in the specialized reports by Gavana \cite{testfunctionindex} and Jamil and Yang \cite{jamil2013}. We selected the test functions commonly referenced in the above two reports to ensure that they have well-defined standards and properties. In total, there are 141 commonly referred test functions. For comparison with other optimization algorithms, we also tested the top 10 hardness optimization functions reported in \cite{testfunctionindex}, seven of which are already included in the 141 aforementioned functions. We added 3 additional functions and resulted in a pool of 144 functions. Of the 144 functions, nine functions are discontinuous (including functions: Corana, Csendes, Damavandi, Helical Valley, Keane, PenHolder, SchmitVetters, Step and Tripod), and one has inconsistent numbers of variables and dimensions (function: Cola). Hence, these functions are excluded from our test function pool. We tested all the remaining 134 functions, which cover the characteristics of differentiability, separability, scalability, and modality.
	
	Each of the 134 tested functions is defined by a formula in which some are dependent on a dimension parameter $n$. By letting $n$ be concrete and different positive integers, we have a pool of more than 134 formulas, each of which we call a sample formula. The testing is performed on the sample formulas. As to the value of $n$, we set the following rule: if the GrS algorithm can find the optimization result in 4000 s, then we consider the sample formula as tested. In our experiments, the highest dimension parameter of all the tested sample formulas is 12. In total, 248 sample formulas were tested by automatic implementation under the same parameter settings. 
	
	In the experiments, we found that 61 out of 248 sample formulas (corresponding to 33 functions) provided by Gavana \cite{testfunctionindex} and Jamil and Yang \cite{jamil2013} had some problems. The problems can be classified into the following four categories.\\
	
	(I) The minimum does not match the function value at the provided minimizer(s). This case can be verified by directly inputting the minimizer(s) into the function.\\
	
	(II) The minimizers provided are incomplete. We found more minimizers and verified them too.\\
	
	(III) Our algorithm found a strictly smaller minimum than the provided value \cite{testfunctionindex,jamil2013}. The new minimum is also verified.\\
	
	(IV) The given formula is incorrect. We found this problem by comparing  Gavana \cite{testfunctionindex} to Adorio and Diliman \cite{Adorio2005}.\\
	
	The tested functions in the experiments are illustrated in the webpage of the Experiment Results \footnote{ https://www.fst.um.edu.mo/personal/lmzhang/experiment\_results\label{web}}, and the problem types are labelled on the problematic functions. The functions for which our algorithm was unsuccessful are marked with “*”.\\

	\noindent{\bf Experimental settings}
	The experiments were performed with MATLAB R2017b on a desktop computer with an Intel Core i9-10920X CPU at a 3.5 GHz clock frequency. In the experiments, for all the 2- and 3-dimensional functions, the side length along each dimension was split into 60 equal segments. For all the functions with more than three dimensions, the side length along each dimension was split into 2 equal segments. The stopping criteria on each $M_i, (i=1,2,...,n)$ was set to $\delta_k M_i \leq 0.001$ or $\delta_k \leq 0.001$ holds for the first time. 
	
	In practice, it is hard to verify whether $M_n$ is an upper bound. To ensure the correctness of the results, when $f_{min}^{M_n}=f_{min}^{M_{n-1}}$ is observed, a few more GrS runs can be applied to ensure optimality. In our experiments, one more run is applied on $M_{n+1}=(2^{n-1}+1)M_1$. This is a trade-off process because the computational cost rapidly increases as the GrS iteration number increases. In summary, if  $f_{min}^{M_{n+1}}=f_{min}^{M_n}=f_{min}^{M_{n-1}}$, then the multiple runs of the GrS algorithm finally stop. In this situation, GrS is ran at least three times in our experiments.\\
	
	\noindent{\bf Results and discussions}	
	Among all the tested formulas, GrS correctly yielded both the minimum and minimizers for 224 sample formulas, corresponding to 119 functions. In the test results of 11 functions (19 sample formulas), the minimums yielded by GrS approximate the ground truths, but the minimizers are incorrect. The analysis shows that this inconsistency is caused by an improperly set termination threshold or an underestimated pseudo-$M$ value. By employing smaller termination thresholds or larger $M$ values, or by allowing computing times greater than 4000 s, the problem of the inconsistency can be solved. The test results of four other functions are different from both the reference provided minimum and minimizers. The analysis of the geometric nature of the function graphs shows that the failures are mostly caused by the first estimation of the pseudo-$M$ in use. This suggests that a larger pseudo-$M$ should be used to increase the effectiveness of the algorithm. With the improved $M$ values, one can successfully find the correct solutions. The success rates of the tested functions and formulas are 88.81\% and 90.32\%, respectively, as listed in Table \ref{t1}. If, like the other studies \cite{testfunctionindex,Liberti2006,Torn1989}, only finding the global minimum is regarded as a success, then our success rate on testing formulas is 98\%.
	\begin{table}[h]
		\caption{Success rates of the GrS algorithm in finding the global optimal solutions.}
		\centering
		\label{t1}
		\begin{tabular}{cccc}
			\hline
			& \textbf{Tested} & \textbf{Successful} & \textbf{Success   rate (\%)} \\ \hline
			Functions & 134             & 119                 & 88.81                        \\
			Formulas  & 248             & 224                 & 90.32                        \\ \hline
		\end{tabular}
	\end{table}
	
	To compare our algorithm with existing ones, in Table \ref{t2}, we illustrate our test results on the top 10 benchmark optimization functions listed in the ``Hardness'' table of Gavana \cite{testfunctionindex}, which provides algorithm comparison among 12 popular global optimization algorithms. The overall success rates in the table were reported to be obtained by running all the available global optimizers against all the test functions for a collection of 100 random starting points and then averaging the successful minimizations across all the optimizers \cite{testfunctionindex}. As presented in Table \ref{t2}, the hardest function for which our algorithm succeeds is XinSheYang03 (No. 4 in the top 10 ``Hardness'' list), and the overall success rate for XinSheYang03 is merely 1.08\% for all the other 12 global optimizers mentioned in Gavana \cite{testfunctionindex}. The non-differentiability of the function makes unavailability of the gradient-dependent methods. This is one of the reasons for the low success rate. This example shows that the gradient-free GrS is an effective algorithm regardless differentiability properties of the objective function. By adding the overall success rate of the top 10 hardest functions in Table 2, we can conclude that the overall success rate of the other 12 popular optimizers on the average of the 10 hardest problems is 3.366\%, while the overall success rate of GrS on the average of the 10 hardest problems is 60\%. 
	
	\begin{table}[h]
		\scriptsize
		\caption{Test results for the top 10 hardness optimization functions listed in Gavana \cite{testfunctionindex}. Test functions with $*$ indicate that there is a problem with the formula or its ground truth.}
		\label{t2}
		\centering
		\begin{tabular}{cccc}
			\hline
			\textbf{\begin{tabular}[c]{@{}c@{}}Hardness \\ ranking\end{tabular}} & \textbf{Test function}  & \textbf{\begin{tabular}[c]{@{}c@{}}Overall success \\rate against \\ground truth (\%)\end{tabular}} & \textbf{\begin{tabular}[c]{@{}c@{}}Result of \\ GrS \end{tabular}} \\ \hline
			1                                                                         & DeVilliersGlasser02                     & 0.00                                                                                               & Failure                          \\ 
			2                                                                         & Damavandi                                & 0.25                                                                                               & Discontinuous                    \\ 
			3                                                                         & CrossLegTable                           & 0.83                                                                                               & Failure                          \\ 
			4                                                                         & XinSheYang03                            & 1.08                                                                                               & Success                          \\ 
			5                                                                         & SineEnvelope$*$                            & 2.17                                                                                               & Success                           \\ 
			6                                                                         & Whitley                                 & 4.92                                                                                               & Success                          \\ 
			7                                                                         & Zimmerman$*$                               & 4.92                                                                                               & Success                      \\ 
			8                                                                         & Griewank                                & 6.08                                                                                               & Failure                          \\ 
			9                                                                         & Trefethen                               & 6.58                                                                                               & Success                          \\ 
			10                                                                        & Bukin06                                  & 6.83                                                                                               & Success                          \\ \hline
		\end{tabular}
	\end{table}
	
\begin{figure}[t]
		\centering
		\includegraphics[scale=0.4]{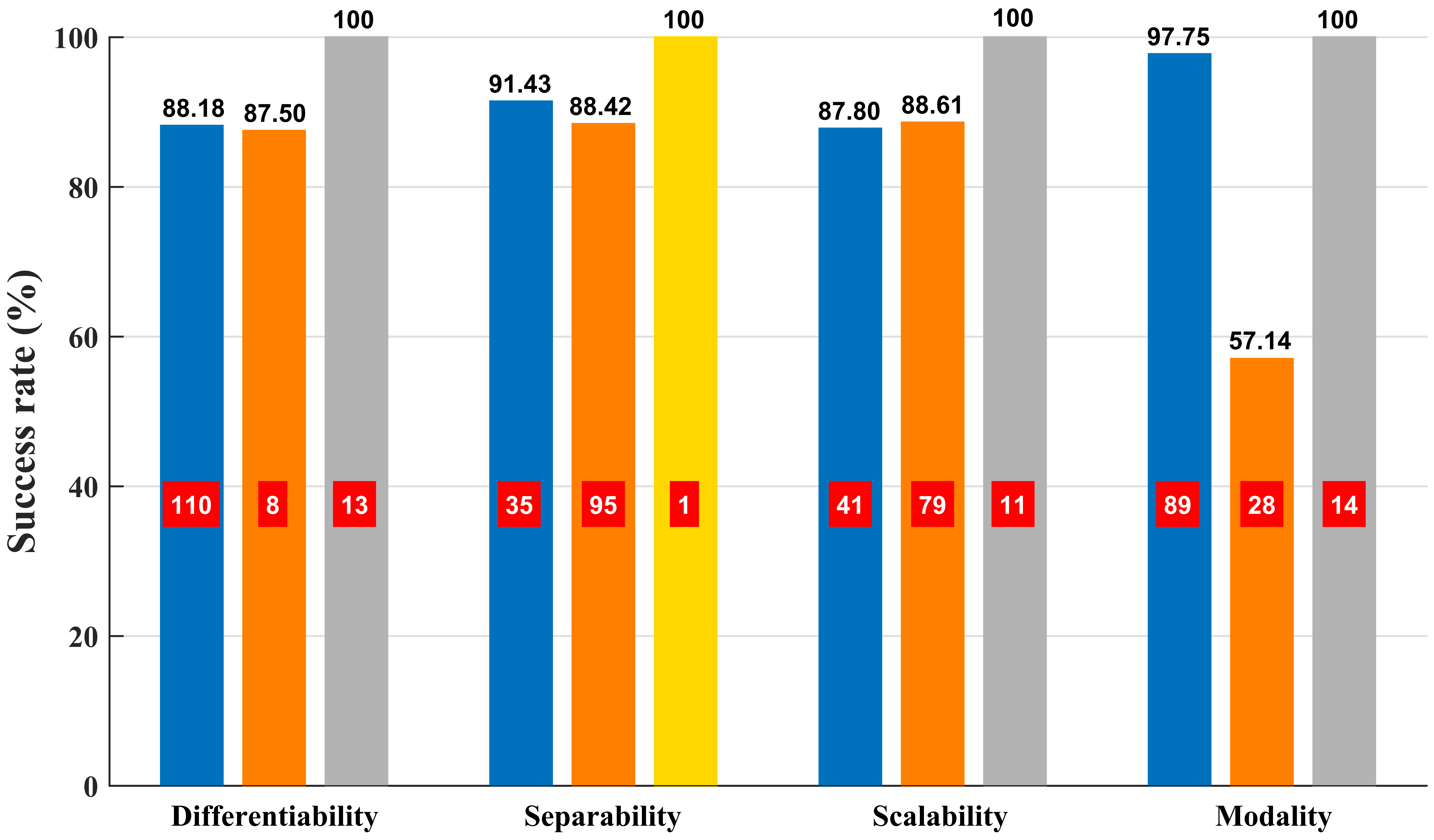}
		\caption{Statistical results on the characteristics of the test functions. The bars with different colors correspond to different types of characteristics. The blue bars represent differentiable, separable, scalable and unimodal types, respectively. The orange ones represent non-differentiable, non-separable, non-scalable and multimodal types. The gray ones are with unmentioned types. The yellow one is the partially-separable type. The white numbers with red background are the corresponding test function numbers.}
		\label{f3}
	\end{figure}

	Apart from the continuity and finite-$M_0$ property, the GrS algorithm does not require any other properties from the objective functions. Our tested functions include non-differentiable, non-separable, scalable and multimodal types, for which global optimal solutions are generally difficult to obtain \cite{jamil2013}. The detailed statistical results on the success rate are presented in Figure \ref{f3}. In total, 131 functions that are commonly referred by \cite{jamil2013} and \cite{testfunctionindex} are listed. The characteristics of the exceptional three functions from Table \ref{t2} are not labeled so that they are not counted. 
	
	Pseudo-$M$ is an important parameter used and estimated along with the algorithm process. In 201 out of 224 successfully tested formulas, $M$ is estimated thrice. Among all the tested sample formulas, the maximal number of estimations for $M$ is 11. The detailed experiment results on all the sample formulas can be found in the webpage of the Experiment Results\textsuperscript{\ref {web}}.
	
\section{Conclusions}

	The proposed GrS method is based on intuitive and simple mathematical formulation and has an easily implementable algorithm. The algorithm is effective and efficient for any finite dimension, offering the complete set of solutions within the prescribed accuracy. If the Lipschitz constant $M_0$ is known, then the algorithm has theoretically no limitation; however, it practically requires a large but relative moderate computer storage capacity.
	
	The real limitation of GrS is related to the fact that, in practice, the Lipschitz constant is unknown. Finding an effective difference-quotient bound, or in other words an effective pseudo-$M$, is practical. This problem is, however, common for all the Lipschitzian methods. $M$-dependence of GrS is analyzed in the study. Significant topics for future studies in this direction would include finding effective pseudo-$M$’s and deeper relationship between the used $M$’s and the GrS.

	\bibliographystyle{siamplain}
	\bibliography{myfile}
\end{document}